\documentclass[11pt]{amsart}

\usepackage[margin=1.25in]{geometry}
\usepackage{tikz}
\usetikzlibrary{shapes.misc}
\usepackage{amsmath,amscd}
\usepackage{amssymb}
\usepackage{amsthm}
\usepackage{oldgerm}
\usepackage[pdftex,hidelinks]{hyperref}
\usepackage{multicol}
\tikzset{cross/.style={cross out, draw=black, minimum size=2*(#1-\pgflinewidth), inner sep=0pt, outer sep=0pt},
cross/.default={1pt}}

\newcommand{\ra}{{\rightarrow}}
\newcommand{\Z}{\mathbb{Z}}

\newcommand{\sF}{\mathcal{F}}

\newcommand{\End}{{\mathbb{E}{{\sf nd}}}}

\newcommand{\BV}{{\sf BV}}
\newcommand{\Lie}{{\sf Lie}}

\newcommand{\mA}{\ensuremath{\mathcal{A}}}
\newcommand{\mO}{\ensuremath{\mathcal{O}}}

\newtheorem{lem}{Lemma}[section]
\newtheorem{thm}{Theorem}[section]

\newtheorem{cor}{Corollary}[section]

\newtheorem{rmk}{Remark}[section]

\title{The Bogomolov-Tian-Todorov theorem of cyclic $A_\infty$-algebras}
\author{Junwu Tu}
\thanks{Institute of Mathematical Sciences, ShanghaiTech University, Shanghai, China, 201210. E-mail: {\texttt{tujw.at.shanghaitech.edu.cn}}}

\begin{document}
\maketitle

\begin{abstract}

Let $A$ be a finite-dimensional smooth unital cyclic $A_\infty$-algebra. Assume furthermore that $A$ satisfies the Hodge-to-de-Rham degeneration property. In this short note, we prove the non-commutative analogue of the Bogomolov-Tian-Todorov theorem: the deformation functor associated with the differential graded Lie algebra of Hochschild cochains of $A$ is smooth. Furthermore, the deformation functor associated with the DGLA of cyclic Hochschild cochains of $A$ is also smooth.

\end{abstract}

\section{The non-commutative Bogomolov-Tian-Todorov Theorem}

Let $X$ be a smooth projective Calabi-Yau variety over a field $\mathbb{C}$. It is a classical result of Bogomolov-Tian-Todorov~\cite{Tia,Tod} that the formal deformation functor associated with the differential graded Lie algebra (DGLA)
\[\mathfrak{g}_X:=\bigoplus_{p,q} \mA^{0,q}(\Lambda^pT_X)\]
of the Dolbeault resolution of holomorphic poly-vector fields is smooth. In order to prove this, the key observation was the existence of a BV operator $ \Delta: \Lambda^*T_X \ra \Lambda^{*-1} T_X$, which ``trivializes" the Lie bracket by the Tian-Todorov identity
\begin{equation}~\label{eq:tian-todorov}
 [\alpha, \beta] = \Delta (\alpha\wedge\beta) - \Delta (\alpha) \wedge \beta - (-1)^{|\alpha|} \alpha\wedge \Delta (\beta).
 \end{equation}
With the above formula, the smoothness of the deformation functor follows easily from the classical $\partial\overline{\partial}$-Lemma in Hodge theory. 

Following Kontsevich-Soibelman~\cite{KS} and Katzarkov-Kontsevich-Pantev~\cite{KKP}, one can formulate the compactness, smoothness, and the Calabi-Yau property purely in terms of the differential graded category of coherent sheaves on $X$. Thus, a natural question is whether the analogues of the Bogomolov-Tian-Todorov's theorem holds for any smooth and proper Calabi-Yau categories. This question might have been a folklore theorem for experts in the field. The purpose of this note is to fill some of the missing details in the literature.

A large class of dg categories of interest is compactly generated by a single object. For this reason, instead of considering formal deformations of dg categories (whatever that means), we shall consider deformations of $A_\infty$-algebras which is also much more tractable. More precisely, we shall let $A$ be an $A_\infty$-algebra over a field $\mathbb{K}$ of characteristic zero. In the geometric context, such an algebra $A$ often appears as 
\[ A= {\sf End}_X^\bullet(E),\]
with $E$ a choice of compact generator of $D^b_{dg}({\sf Coh}(X))$.  If $A$ is a chain complex, its suspension is denoted by $sA$ with $(sA)_n:= A_{n-1}$. A {\em cyclic structure} on $A$ of dimension $d$ is a symmetric non-degenerate pairing $\langle -,-\rangle:(sA)^{\otimes 2}\to \mathbb{K}$ of degree $d-2$ that satisfies 
\begin{equation}\label{eq:cyclic-structure}
\langle m_k(a_0, \ldots, a_{k-1}), a_k\rangle=(-1)^{|a_0|'(|a_1|'+\ldots+|a_k|')}\langle m_k(a_1, \ldots, a_{k}), a_0\rangle,
\end{equation}
where $|a_i|'=|a_i|-1$. We refer to~\cite{Pol} for a detailed explanation how in the geometric context, a holomorphic volume form (given by the Calabi-Yau condition) gives rise to a cyclic structure on ${\sf End}_X^\bullet(E)$.

To state the non-commutative version of the Bogomolov-Tian-Todorov theorem precisely, we first fix some notations and conventions. Throughout the paper, we use the homological degree of chain complexes. For a unital $A_\infty$-algebra $A$, denote by $C_{mod}^{-*}(A)$ ($C^{mod}_*(A)$) its modified Hochschild cochain complex (chain complex respectively), defined by Kontsevich-Soibelman in~\cite[Section 11.3]{KS}. The minus sign on degree is due to that we use homological degree. Let $A$ be a cyclic unital $A_\infty$-algebra, denote by $C^\lambda(A)\subset C^{-*}_{mod}(A)$ the sub-complex consisting of Hochschild cochains that are cyclic (i.e. satisfying Equation~\eqref{eq:cyclic-structure} with respect to the pairing on $A$. 

\medskip

\begin{thm}~\label{thm:main}
Let $A$ be a $\Z/2\Z$-graded, finite-dimensional smooth cyclic $A_\infty$-algebra over a field $\mathbb{K}$ of characteristic zero. Assume furthermore that $A$ satisfies the Hodge-to-de-Rham degeneration property. Then we have
\begin{itemize}
\item[(A.)] The deformation functor ${{\sf Def}}_{sC^{-*}_{mod}(A)}$ associated with the DGLA $sC^{-*}_{mod}(A)$ is smooth.
\item[(B.)] The deformation functor ${{\sf Def}}_{sC^\lambda(A)}$ associated with the DGLA $sC^\lambda(A)$ is also smooth.
\item[(C.)] The natural transformation $f: {{\sf Def}}_{sC^\lambda(A)}\ra {{\sf Def}}_{sC^{-*}_{mod}(A)}$ associated with the canonical inclusion map $sC^\lambda(A)\hookrightarrow sC^{-*}_{mod}(A)$ is smooth. In particular, every deformation of the $A_\infty$ structure of $A$ lifts to a deformation of the {{\em cyclic}} $A_\infty$ structure of $A$.
\end{itemize}
\end{thm}

\begin{rmk}
The assumption of the Hodge-to-de-Rham degeneration property automatically holds for any $\Z$-graded smooth and proper $A_\infty$-algebra by Kaledin~\cite{Kal}. In the general $\Z/2\Z$-graded case, this remains an open conjecture by~\cite{KKP,KS}. 
\end{rmk}
\begin{rmk}
 Part (A.) of the theorem in the $\Z$-graded case was proved by Isamu Iwanari \cite{Iwa} using derived stacks and period mappings. In contrast, our approach utilizes the $BV_\infty$ machinery of the framed little disks operad to provide a somewhat more explicit proof. We also point out that while in this paper we focus on single-object algebras, the results extend to $A_\infty$-categories with minor changes, see~\cite{AT} for the appropriate categorical constructions.
\end{rmk}

\section{$\BV_\infty$-algebra structure on $C^{-*}_{mod}(A)$}

To prove the above Theorem~\ref{thm:main}, one follows the same idea as in the proof of Bogomolov-Tian-Todorov Theorem. However, the key identity~\eqref{eq:tian-todorov} fails to hold. It only holds up to homotopy. Also, the cup product on $C^{-*}_{mod}(A)$ is only commutative up to homotopy. Thus it is natural to work with a homotopy version of the underlying algebraic structures. 

In this section, we first exhibit a homotopy BV algebra structure, or $\BV_\infty$-algebra structure on $C^{-*}_{mod}(A)$. The definition of $\BV_\infty$-algebras used in this paper is from the article~\cite{GCTV}. In fact, it was argued in {{\sl Loc. Cit.}} that combining a TCFT structure defined by~\cite{Cos}~\cite{KS} and the formality of the operad $\BV$, one easily deduces the existence of a $\BV_\infty$-algebra structure on $C^{-*}_{mod}(A)$ with $A$ as in Theorem~\ref{thm:main}. However, to make such a structure useful in order to deduce Theorem~\ref{thm:main}, one needs to say a bit more about this $\BV_\infty$ structure. For example, its underlying $\Lie_\infty$ algebra is in fact given by the differential graded Lie algebra $\big( C^{-*}_{mod}(A), \delta, [-,-]_G\big)$. For this purpose, we need to use a construction of Tamarkin in his proof of the Deligne's conjecture~\cite{Tam}. We introduce the following notations:
\begin{itemize}
\item $\Lie$ --- The Lie operad.
\item $\Lie_\infty$ --- The homotopy Lie operad.
\item $E_2$ --- The operad whose representation gives Gerstenhaber algebras.
\item $G_\infty$ --- The homotopy $E_2$ operad.
\item $\BV$ --- The operad whose representation gives BV algebras.
\item $\BV_\infty$ --- The homotopy BV operad.
\item $B_\infty$ --- The brace operad~\cite[Section 5.2]{GJ}.
\item $\sF$ --- The operad defined by Tamarkin in~\cite[Section 6]{Tam}.
\item $C^{{\sf comb}}_*(FD)$ --- The operad of black-and-white ribbon trees defined by Kontsevich-Soibelman~\cite[Section 11.6]{KS}, see also Wahl-Westerland~\cite[Section 2]{WW}. This operad gives a combinatorial model for the framed little disk operad.
\item ${{\mathcal{X}}}$ --- A cofibrant replacement of $ C^{{\sf comb}}_*(FD)$.
\end{itemize}

For an operad $\mO$, denote by $\mO\left\{1\right\}$ its shifted version so that an $\mO\left\{1\right\}$-algebra structure on a chain complex $A$ is equivalent to an $\mO$-algebra structure on $sA$. The endomorphism operad of a chain complex is denoted by $\End(-)$.

The starting point to construct a $\BV_\infty$ structure on $C^{-*}_{mod}(A)$ is that the operad $C^{{\sf comb}}_*(FD)$ naturally acts on $C^{-*}_{mod}(A)$:
\begin{equation}~\label{eq:action} C^{{\sf comb}}_*(FD) \ra \End\big( C^{-*}_{mod}(A) \big).\end{equation}
We refer to Kontsevich-Soibelman~\cite{KS} and Wahl-Westerland~\cite{WW} for details of this action.  Here we illustrate this action with a few examples. Indeed, the following black-and-white ribbon tree
\[ \begin{tikzpicture}[baseline={([yshift=-1ex]current bounding
		box.center)},scale=0.2]
\draw (0,0) node[cross=5pt,label=above:{}] {};
\draw[thick]  (0,0) to (2,0);
\draw (2.5,0) circle (.5);
\draw [ultra thick] (3,0) to (5,0);
\end{tikzpicture}\]
inside $C^{{\sf comb}}_*(FD)(1)$ gives the pull-back of the Connes operator $$B: C_{*}^{mod}(A) \ra C_{*}^{mod}(A)$$ under the isomorphism $C^{-*}_{mod}(A)\cong \big( C_{*}^{mod}(A) \big)^\vee$. Denote this operator by $\Delta: C^{-*}_{mod}(A)\ra C^{-*}_{mod}(A)$. There are also binary operators associated with trees in $C^{{\sf comb}}_*(FD)(2)$. For example, consider the following two black-and-white ribbon trees:
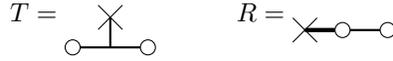
\[  T= \begin{tikzpicture}[baseline={([yshift=-1ex]current bounding
		box.center)},scale=0.2]
\draw (2,2) node[cross=5pt,label=above:{}] {};
\draw[thick]  (0,0) to (4,0);
\draw[thick] (2,2) to (2,0);
\draw (-.5,0) circle (.5);
\draw (4.5,0) circle (.5);
\end{tikzpicture}\;\;\;\;\;\;\;\;\;\;\;R=\begin{tikzpicture}[baseline={([yshift=-1ex]current bounding
		box.center)},scale=0.2]
\draw (0,0) node[cross=5pt,label=above:{}] {};
\draw[ultra thick]  (0,0) to (2,0);
\draw (2.5,0) circle (.5);
\draw [thick] (3,0) to (5,0);
\draw (5.5,0) circle (.5);
\end{tikzpicture}\]
The graph $T$ gives the familiar cup product on $C^{-*}_{mod}(A)$, while the graph $R$ gives the first brace operator $\circ$ on $C^{-*}_{mod}(A)$ whose commutator (after shift) is the Gerstenhaber Lie bracket $[-,-]_G$.

\medskip
\begin{lem}\label{lem:brace}
There is a morphism of operads $h: B_\infty \ra C^{{\sf comb}}_*(FD)$ defined by
\begin{align*}
 h(m_k)&:= T_k, k\geq 2 \\
 h(m_{1,k} )&:= W_k, k\geq 1,\\
 h(m_{p,k}) &:= 0, p\geq 2.
 \end{align*}
Here $T_k$ and $W_k$ are given by the following black-and-white ribbon trees:
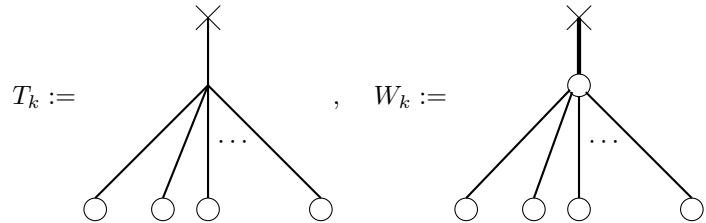
\[T_k := \begin{tikzpicture}[baseline={(current bounding box.center)},scale=0.3]
\draw (0,8) node[cross=5pt,label=above:{}] {};
\draw (-5,-.5) circle (.5);
\draw (-2,-.5) circle (.5);
\draw (0,-.5) circle (.5);
\draw (5,-.5) circle (.5);
\node at (1.2,2.5) {$\cdots$};
\draw [thick] (0,5) to (-5,0);
\draw [thick] (0,5) to (-2,0);
\draw [thick] (0,5) to (0,0);
\draw [thick] (0,5) to (5,0);
\draw [thick] (0,8) to (0,5);
\end{tikzpicture}, \;\;\;\;W_k:=\begin{tikzpicture}[baseline={(current bounding box.center)},scale=0.3]
\draw (0,8) node[cross=5pt,label=above:{}] {};
\draw (-5,-.5) circle (.5);
\draw (-2,-.5) circle (.5);
\draw (0,-.5) circle (.5);
\draw (5,-.5) circle (.5);
\draw (0,5) circle (.5);
\node at (1.2,2.5) {$\cdots$};
\draw [thick] (-.4,4.8) to (-5,0);
\draw [thick] (-.3,4.7) to (-2,0);
\draw [thick] (0,4.5) to (0,0);
\draw [thick] (.3,4.7) to (5,0);
\draw [ultra thick] (0,8) to (0,5.5);
\end{tikzpicture} \]
\end{lem}

\proof This is a straight-forward check.\qed

\medskip
In~\cite{Tam}, Tamarkin constructed morphisms $s: G_\infty \ra \sF$ and $t: \sF \ra B_\infty$. 

\medskip
\begin{lem}~\label{lem:lift}
There exists a commutative diagram:
\[\begin{tikzpicture}
\node at (0,0) {$\sF$};
\node at (1.2, .2) {$t$};
\draw [thick,->] (.3,0) to (2,0);
\node at (2.4,0) {$B_\infty$};
\draw [thick,->] (2.4, -.2) to (2.4, -2.2);
\node at (2.6,-1.2) {$h$};
\node at (3.2, -2.5) {$C^{{\sf comb}}_*(FD)$};
\node at (0,-2.5) {${\mathcal{X}}$};
\draw [thick,dashed,->] (0,-.3) to (0, -2.2);
\draw [thick,->] (.3, -2.5) to (2.2, -2.5);
\end{tikzpicture}\]
\end{lem}

\proof This follows from the lifting property since ${\mathcal{X}} \ra C^{{\sf comb}}_*(FD)$ is a trivial fibration, while $\sF$ by construction is cofibrant. \qed

\medskip
The framed little disk operad is known to be formal with cohomology the $\BV$ operad, which implies that $\BV\cong C^{{\sf comb}}_*(FD)$ in the homotopy category of differential graded operads. Since the operad ${\mathcal{X}}$ is cofibrant, and $\BV$ is fibrant (as any dg operad is fibrant), we obtain a morphism ${\mathcal{X}}\ra \BV$ such that the roof diagram
\[\begin{tikzpicture}
\node at (0,0) {$C^{{\sf comb}}_*(FD)$};
\node at (4,0) {$\BV$};
\node at (2, 2) {$\mathcal{X}$};
\draw [thick,->] (1.8,1.8) to (0,.3);
\draw [thick,->] (2.2,1.8) to (4,.3);
\end{tikzpicture}\]
represents an isomorphism in the homotopy category of differential graded operads. 

\medskip
\begin{lem}
The following diagram is commutative:
\[\begin{tikzpicture}
\node at (-.5,0) {$\Lie_\infty\left\{1\right\}$};
\draw [thick,->] (.3,0) to (.8,0);
\node at (1.2,0) {$G_\infty$};
\draw [thick,->] (1.5,0) to (2.2,0);
\node at (2.5,0) {$\sF$};
\draw [thick,->] (2.8,0) to (3.5,0);
\node at (3.7,0) {$\mathcal{X}$};
\draw [thick,->] (3.7, -.2) to (3.7, -2.2);
\node at (3.9, -2.5) {$\BV$};
\node at (0,-2.5) {$\BV_\infty$};
\draw [thick,->] (0,-.3) to (0, -2.2);
\draw [thick,->] (.3, -2.5) to (3.5, -2.5);
\end{tikzpicture}\]
\end{lem}

\proof It is clear that the left composition factors as
\[ \Lie_\infty\left\{1\right\} \ra \Lie\left\{1\right\} \ra BV.\]
For the right composition, consider the following composition
\[ \Lie_\infty\left\{1\right\} \ra \sF \ra \mathcal{X} \ra C_*^{{\sf comb}}(FD).\]
By Lemma~\ref{lem:lift} above, it is equal to
\[ \Lie_\infty\left\{1\right\} \ra \sF \ra B_\infty \ra C_*^{{\sf comb}}(FD),\]
which by~\cite{Tam} can be factored as
\[ \Lie_\infty\left\{1\right\} \ra \Lie\left\{1\right\} \ra B_\infty \ra C_*^{{\sf comb}}(FD).\]
This shows that both compositions vanish on generators $l_k$ of $\Lie_\infty\left\{1\right\}$ with $k\geq 3$. For $k=2$, this is a direct check by definition.\qed

\medskip
By definition of ${\mathcal{X}}$, the right vertical map in the above diagram is a trivial fibration since $\BV$ is minimal. The left vertical map is a cofibration. Thus by the lifting property, we obtain a map $\BV_\infty \ra {\mathcal{X}}$:
\[\begin{tikzpicture}
\node at (-.5,0) {$\Lie_\infty\left\{1\right\}$};
\draw [thick,->] (.3,0) to (.8,0);
\node at (1.2,0) {$G_\infty$};
\draw [thick,->] (1.5,0) to (2.2,0);
\node at (2.5,0) {$\sF$};
\draw [thick,->] (2.8,0) to (3.5,0);
\node at (3.7,0) {$\mathcal{X}$};
\draw [thick,->] (3.7, -.2) to (3.7, -2.2);
\node at (3.9, -2.5) {$\BV$};
\node at (0,-2.5) {$\BV_\infty$};
\draw [thick,->] (0,-.3) to (0, -2.2);
\draw [thick,->] (.3, -2.5) to (3.5, -2.5);
\draw [thick,dashed,->] (.3,-2.2) to (3.5, -.3);
\end{tikzpicture}\]
We define a $\BV_\infty$-algebra structure on $C^{-*}_{mod}(A)$ via the composition:
\begin{equation}~\label{eq:bv}
\rho:  \BV_\infty \ra \mathcal{X} \ra C^{{\sf comb}}_*(FD) \ra \End\big( C^{-*}_{mod}(A)\big).
 \end{equation}

\medskip
\begin{cor}~\label{cor:lie}
Let $C^{-*}_{mod}(A)$ be endowed with the $BV_\infty$-algebra structure defined by Equation~\eqref{eq:bv}. Then its underlying $\Lie_\infty$ structure on the suspension of $C^{-*}_{mod}(A)$ is given by the differential graded Lie algebra $\big( s C^{-*}_{mod}(A), \delta, [-,-]_G\big)$.
\end{cor}

\proof This follows from Lemma~\ref{lem:lift} and Tamarkin's commutative diagram:
\[\begin{tikzpicture}
\node at (-.5,0) {$\Lie_\infty\left\{1\right\}$};
\draw [thick,->] (.3,0) to (2,0);
\node at (2.4,0) {$G_\infty$};
\draw [thick,->] (2.8,0) to (4.5,0);
\node at (4.8,0) {$\sF$};
\draw [thick,->] (4.8, -.2) to (4.8, -2.2);
\node at (4.8, -2.5) {$B_\infty$};
\node at (-.5,-2.5) {$\Lie\left\{1\right\}$};
\draw [thick,->] (0,-.3) to (0, -2.2);
\draw [thick,->] (.3, -2.5) to (4.5, -2.5);
\end{tikzpicture}\]
\qed

Consider a subset of generating operators of a $\BV_\infty$ structure given by $$\left\{ l_k^d | k\geq 1, d\geq 0\right\}$$ which label a basis of the convolution between the Koszul dual cooperad of $\Lie$ and that of the operad generated by the circle operator $\Delta$. We denote the sub-operad generated by $l_k^d$ in $\BV_\infty$ by ${{\sf q}}\Lie_\infty\left\{1\right\}$. The notation is because a ${{\sf q}}\Lie_\infty\left\{1\right\}$-algebra structure on a chain complex $V$ is equivalent to an $\Lie_\infty$ structure on $sV[[u]]$ (with $u$ a degree $2$ formal variable), which may be thought of as a ``quantum" $\Lie_\infty\left\{1\right\}$ structure. 

\medskip
\begin{cor}~\label{cor:qlie}
The induced ${{\sf q}}\Lie_\infty\left\{1\right\}$ structure on $C^{-*}_{mod}(A)$ is of the form
\begin{align}
& l_1^0=\delta, \; l_2^0=[-,-]_G, \; l_k^0=0 \; \forall \; k\geq 3, \\
& l_1^1=\Delta, \; l_k^1(\alpha_1,\cdots,\alpha_k)=0, \forall \alpha_1,\cdots,\alpha_k\in {{\sf Ker}} \Delta, \;\; \forall k\geq 2,\\
& l_k^d=0, \; \forall \; d\geq 2.
\end{align}
\end{cor}

\proof Property $(3)$ is proven in the previous Corollary~\ref{cor:lie}. To prove Property $(4)$ and $(5)$, observe that the degree of the operator $l_k^d$ is equal to $2d+2k-3$, see~\cite{GCTV}. But, a top dimensional chain in the chain complex $C^{{\sf comb}}_*(FD)(k)$ has degree equal to $2k-1$, which implies the vanishing in $(5)$. The operator $l_1^1=\Delta$, since the black-and-white graph $$ \begin{tikzpicture}[scale=0.2]
\draw (0,0) node[cross=5pt,label=above:{}] {};
\draw[thick]  (0,0) to (2,0);
\draw (2.5,0) circle (.5);
\draw [ultra thick] (3,0) to (5,0);
\end{tikzpicture}$$ is the unique degree one graph representing the fundamental class in the homology group $H_1\big( C^{{\sf comb}}_*(FD)(1) \big)\cong H_1(\mathbb{S}^1)$. For $k\geq  2$, since the degree of the operator $l_k^1$ is $2k-1$ which is top dimensional in the chain complex $C_*^{{\sf comb}}(FD)(k)$. Thus $l_k^1$ is a linear combination of operators coming from black-and-white ribbon trees with $k$ white vertices, and of degree $2k-1$. One can show that any such tree must contain at least one white vertex of the form:
\[\begin{tikzpicture}[scale=0.2]
\node at (-.4,0) {$\cdots$};
\draw[thick]  (.5,0) to (2,0);
\draw (2.5,0) circle (.5);
\draw [ultra thick] (3,0) to (5,0);
\end{tikzpicture}\]
Since otherwise, the maximal degree on such a tree would be 
\[ \sum_{i=1}^k ({\sf val}(w_i) -1)= \sum_{i=1}^k {\sf val}(w_i) -k = 1+2(k-1) -k = k-1 < 2k-1,\]
which would be a contradiction. If we input any $\alpha\in {{\sf Ker}}\Delta$ at the type of white vertex shown in the figure above, the result gives zero by definition of the action of ribbon trees on $C^{-*}_{mod}(A)$ (by Equation~\eqref{eq:action}).\qed

\section{Proof of Theorem~\ref{thm:main}}

As explained in the previous section, via the composition in Equation~\eqref{eq:bv}, we have a $\BV_\infty$ structure on $C^{-*}_{mod}(A)$. Thus, we may form its bar-cobar resolution~\cite{LV}:
\[ \Omega B C^{-*}_{mod}(A) \cong C^{-*}_{mod}(A),\]
which yields a differential graded BV algebra homotopy equivalent to $C^{-*}_{mod}(A)$. At this point, we use the following theorem due to Katzarkov-Kontsevich-Pantev~\cite{KKP} and Terilla~\cite{Ter}.

\begin{thm}
Let $S$ be a differential graded BV algebra. Assume that the spectral sequence associated with the complex $\big( S[[u]], d+u\Delta\big)$ endowed with the $u$-filtration is degenerate at the first page. Then both the DGLAs $\big( S, d, [-,-]\big)$ and $\big( S[[u]], d+u\Delta, [-,-]\big)$ are homotopy abelian.
\end{thm}

Since the degeneration of the spectral sequence is a homotopy invariant property (see~\cite{DSV}), we may use the above theorem to deduce the homotopy abelian property of the DGLA $sC^{-*}_{mod}(A)$. Note that here it is essential that the $\BV_\infty$ structure on $C^{-*}_{mod}(A)$ extends the DGLA structure of $sC^{-*}_{mod}(A)$ by Corollary~\ref{cor:lie}.

Similarly, we may also restrict the $\BV_\infty$ structure on $C^{-*}_{mod}(A)$ to the sub-operad ${{\sf q}}\Lie_\infty\left\{1\right\}$. Then the theorem above implies that the following $\Lie_\infty$-algebra
\[ \big( sC^{-*}_{mod}(A)[[u]], \;\delta+u\Delta,\; [-,-]_G+u\cdot l^1_2,\; u\cdot l^1_k \; (k\geq 3) \big)\]
is homotopy abelian. Here the structure maps $l^1_k$ are as in Corollary~\ref{cor:qlie}. The following lemma then finishes the proof of $(B.)$ in Theorem~\ref{thm:main}, using the homotopy invariance of deformation functors. 

\begin{lem}
The canonical inclusion map 
\[ \iota: sC^\lambda(A) \hookrightarrow sC^{-*}_{mod}(A)[[u]]\]
is a quasi-isomorphism of $\Lie_\infty$-algebras.
\end{lem}

\proof This inclusion is a quasi-isomorphism is a classical result, see for example~\cite[Section 2.4]{Lod}.  The fact that the Gerstenhaber bracket of two cyclic cochains remains a cyclic cochain was proved in~\cite{Pen}. By Corollary~\ref{cor:qlie}, the higher brackets $l^1_k$ vanish on elements inside ${{\sf Ker}}\Delta$ while we certainly have that ${{\sf Im}}\,\iota\subset {{\sf Ker}}\Delta$. This shows that
\[ l^1_k(\iota \alpha_1,\cdots,\iota\alpha_k)=0, \;\; \forall \alpha_1,\cdots,\alpha_k\in sC^\lambda(A),\]
which implies that $\iota$ is a morphism of $\Lie_\infty$-algebras.\qed

\medskip
Denote by $\pi: sC^{-*}_{mod}(A)[[u]]\ra sC^{-*}_{mod}(A)$ the projection map defined by setting $u=0$. Part $(C.)$ of Theorem~\ref{thm:main} easily follows from $(A.)$ and $(B.)$. Indeed, since both functors are smooth, it suffices to check that the inclusion map
\[ \pi\circ \iota: sC^\lambda(A) \hookrightarrow sC^{-*}_{mod}(A)\]
induces a surjection on the tangent spaces of the associated deformation functors. This is clear as $\iota$ is a quasi-isomorphism, and $\pi$ is a surjective on cohomology by the Hodge-to-de-Rham degeneration assumption.

\end{document}